\documentclass[12pt]{amsart}
\usepackage[all]{xy} % for complicated commutative diagrams

\def\ra{\rightarrow}
\def\e{\kern 0.08em}
\def\ng{\kern -0.03em}
\def\be{\kern -0.015em}
\def\g{\varGamma}
\def\kb{\overline{k}}
\def\cb{\overline{C}}

\def\xb{\overline{X}}

\def\kc{\kb\!\phantom{.}^{*}}

\newtheorem{theorem}{Main Theorem\!\!}

\newtheorem{lemma}{Lemma}[section]
\newtheorem{teorema}[lemma]{Theorem}
\newtheorem{corollary}[lemma]{Corollary}
\newtheorem{proposition}[lemma]{Proposition}
\theoremstyle{definition}

\theoremstyle{remark}

\begin{document}

\title[Algebraic cycles on Severi-Brauer schemes]
{Algebraic cycles on Severi-Brauer schemes of prime degree over a
curve}


\author{Cristian D. Gonz\'alez-Avil\'es}
\address{Departamento de Matem\'aticas, Universidad Andr\'es Bello,
Chile} \email{cristiangonzalez@unab.cl}


\thanks{The author is partially supported by Fondecyt grant
1061209 and Universidad Andr\'es Bello grant DI-29-05/R}

\maketitle

\date{January 3rd, 2007}

\begin{abstract} Let $k$ be a perfect field and let $p$ be a prime
number different from the characteristic of $k$. Let $C$ be a
smooth, projective and geometrically integral $k$-curve and let $X$
be a Severi-Brauer $C$-scheme of relative dimension $p-1$ . In this
paper we show that $CH^{d}(X)_{{\rm{tors}}}$ contains a subgroup
isomorphic to $CH_{0}(X/C)$ for every $d$ in the range $2\leq d\leq
p$. We deduce that, if $k$ is a number field, then $CH^{d}(X)$ is
finitely generated for every $d$ in the indicated range.

\end{abstract}

\section{Introduction.}
Let $k$ be a perfect field with algebraic closure $\kb$. Very little
is known about algebraic cycles on algebraic $k$-varieties,
especially in codimension greater than 2 or dimension greater than
zero. Let $p$ be a prime number different from the characteristic of
$k$ and let $C$ be a smooth, projective and geometrically integral
$k$-curve. In this paper we study a certain subgroup of
$CH^{d}(X)_{{\rm{tors}}}$ for a Severi-Brauer $C$-scheme $q\colon
X\ra C$ of relative dimension $p-1$ and any integer $d$ such that
$2\leq d\leq p$. Let
$$
CH_{0}(X/C)={\rm{Ker}}\!\left[\e CH_{0}(X)\overset{q_{
*}}{\longrightarrow} CH_{0}(C)\e\right]
$$
and let $\pi^{*}\colon CH^{d}(X)\ra CH^{d}\!\left(\e\xb\e\right)$ be
induced by the extension-of-scalars map $\xb\ra X$, where
$\xb=X\otimes_{k}\kb$. Then the following holds.
\begin{theorem} For any $d$ as above, there exists a canonical
isomorphism
$$
{\rm{Ker}}\!\left[\e CH^{d}(X)\overset{\pi^{*}}{\longrightarrow}
CH^{d}\!\left(\e\xb\e\right) \e\right]\simeq CH_{0}(X/C).
$$
Consequently, if $k$ is a number field, then $CH^{d}(X)$ is finitely
generated.
\end{theorem}

\section*{Acknowledgement.}
I thank B.Kahn for some helpful comments.

\section{Preliminaries.}
Let $k$ be a perfect field, fix an algebraic closure $\kb$ of $k$
and let $\g={\rm{Gal}}\!\left(\e\kb/k\right)$. Now let $C$ be a
smooth, projective and geometrically integral $k$-curve and let $X$
be a Severi-Brauer scheme over $C$ [4, \S 8]. There exists a proper
and flat $k$-morphism $q\colon X\ra C$ all of whose fibers are
Severi-Brauer varieties of dimension $m-1$ ($m\geq 1$) over the
appropriate residue field [loc.cit.]. We will write $X_{\eta}$ for
the generic fiber $X\times_{ C}{\rm{Spec}}\e k(C)$ of $q$ and $A$
for the central simple $k(C)$-algebra associated to $X_{\eta}$. We
define
$$
CH_{0}(X/C)={\rm{Ker}}\!\left[\e CH_{0}(X)\overset{q_{
*}}{\longrightarrow} CH_{0}(C)\e\right].
$$
Now let $C_{0}$ be the set of closed points of $C$. The group of
{\it divisorial norms} of $X/C$ (cf. [6]) is the group
$$
k(C)^{*}_{{\rm{dn}}}=\left\{f\in k(C)^{*}\colon\forall\e y\in
C_{0},{\rm{ord}}_{\e y}\be(f)\in (q_{\e
y}\be)_{*}\!\left(CH_{0}(X_{y})\right)\right\}
$$
where, for each $y\in C_{0}$, $q_{\e y}\colon X_{y}\ra{\rm{Spec}}\e
k(y)$ is the structural morphism of the fiber $X_{y}$. This group is
closely related to $CH_{0}(X/C)$ (see [2, Proposition 3.1]). Indeed,
there exists a canonical isomorphism
$$
CH_{0}(X/C)\simeq k(C)^{*}_{{\rm{dn}}}/k^{*}{\rm{Nrd}}\e A^{*}.
$$
Now fix an integer $d$ such that $1\leq d\leq m$ and let
$$
CH^{d}(X)^{\prime}={\rm{Ker}}\!\left[\e
CH^{d}(X)\overset{\pi^{*}}{\longrightarrow}
CH^{d}\!\left(\e\xb\e\right)^{\g}\e\right],
$$
where $\pi\colon\xb\ra X$ is the canonical map. A simple transfer
argument shows that $CH^{d}(X)^{\prime}$ is a subgroup of
$CH^{d}(X)_{{\rm{tors}}}$. Now, since $\xb\ra\cb$ has a section,
$\xb$ is a projective bundle over $\cb$. Thus there exists an
isomorphism
$$
CH^{d}\!\left(\e\xb\e\right)\simeq \Bbb Z\oplus
CH_{0}\!\left(\e\cb\e\right).
$$
(see [3, Theorem 3.3(b), p.34]). Therefore, if $J_{C}(k)$ is
finitely generated, where $J_{C}$ is the Jacobian variety of $C$
(e.g., $k$ is a number field or $C={\Bbb P}^{1}_{k}$), then
$CH^{d}(X)$ is finitely generated if and only if
$CH^{d}(X)^{\prime}$ is finite.

\section{The general method.}
Let $C$ be as above and let $X$ be any smooth, projective and
geometrically integral $k$-variety such that there exists a proper
and flat morphism $q\colon X\ra C$ whose generic fiber $X_{\eta}$ is
geometrically integral. We have an exact sequence [7]
\begin{equation}
H^{d-1}(X_{\eta},{\mathcal
K}_{d})\overset{\delta}{\longrightarrow}\bigoplus_{y\in
C_{0}}CH^{d-1}(X_{y}) \ra CH^{d}(X)\overset{j^{*}}\to
CH^{d}(X_{\eta})\ra 0,
\end{equation}
where $j\colon X_{\eta}\ra X$ is the natural map and the map which
we have labeled $\delta$ will play a role later when $k=\kb$. A
similar exact sequence exists over $\kb$, and we have two natural
exact commutative diagrams:
\[
\xymatrix{0\ar[r]&{\rm{Ker}}\,j^{*}\ar[r]\ar[d]&
CH^{d}(X)\ar[r]\ar[d]&CH^{d}(X_{\eta})\ar[d]
\ar[r]&0\\
0\ar[r]&\left({\rm{Ker}}\,\overline{\jmath}^{\,
*}\right)^{\g}\ar[r]& CH^{d}\!\left(\e\xb\e\right)^{\g}\ar[r] &
CH^{d}\!\left(\e\xb_{\overline{\eta}}\e\right)^{\g}& &
\\
}
\]
and

\begin{equation}
\xymatrix{ 0\ar[r]&\frac{H^{d-1}(X_{\eta},\e {\mathcal
K}_{d})}{j^{*}H^{d-1}(X,\e {\mathcal K}_{d})}
\ar[r]\ar[d]&\displaystyle\bigoplus_{y\in
C_{0}}CH^{d-1}\!\left(X_{y}\right)
\ar[d]^{\oplus_{\e\overline{y}\e\mid y}\pi_{\overline{y}}^{*}}
\ar@{->>}[r] & {\rm{Ker}}\,j^{*}\ar[d]\\
0\ar[r]& \left(\frac{H^{d-1}\!\left(\e\xb_{\overline{\eta}},\e
{\mathcal
K}_{d}\right)}{\overline{\jmath}^{\,*}H^{d-1}\!\left(\e\xb,\e
{\mathcal K}_{d}\right)}\right)^{\g} \ar[r]^{{\kern
-1em}\overline{\delta}}&\displaystyle\bigoplus_{\overline{y}\e\mid
y}CH^{d-1}\!\left(\e \xb_{\overline{y}}\right)^{\g_{\ng y}}{\kern
-.5em}\ar[r] &\left({\rm{Ker}}\,\overline{\jmath}^{\,*}\right)^{\g}
\\
}
\end{equation}
where, for each $y\in C_{0}$, we have fixed a closed point
$\overline{y}$ of $\cb$ lying above $y$ and written $\g_{\ng
y}={\rm{Gal}}\left(\e \kb/k(y)\right)$. Set
$$
CH^{d}\!\left(X_{\eta}\right)^{\e\prime}={\rm{Ker}}\!\left[\e
CH^{d}\!\left(X_{\eta}\right)\ra
CH^{d}\!\left(\e\xb_{\overline{\eta}}\e\right)^{\g}\e\right]
$$
and, for each $y\in C_{0}$,
$$
CH^{d-1}\!\left(X_{y}\right)^{\e\prime}={\rm{Ker}}\!\left[\e
CH^{d-1}\!\left(X_{y}\right)\overset{\pi_{\overline{y}}^{*}}
{\longrightarrow} CH^{d-1}\!\left(\e
\xb_{\overline{y}}\right)^{\g_{\ng y}}\e\right].
$$
Now define
\begin{equation}
E\!\left(\e\xb/\e\cb\,\right)={\rm{Coker}}\!\left[\e\frac{H^{d-1}\!\left(X_{\eta},\e
{\mathcal K}_{d}\right)}{j^{*} H^{d-1}\!\left(X,\e {\mathcal
K}_{d}\right)}\longrightarrow \left(\frac{H^{d-1}\!\left(\e
\xb_{\overline{\eta}},\e {\mathcal
K}_{d}\right)}{\overline{\jmath}^{\e *} H^{d-1}\!\left(\e\xb,\e
{\mathcal K}_{d}\right)}\right)^{\! \g}\,\,\right].
\end{equation}
Then, applying the snake lemma to the preceding diagrams, we
obtain\footnote{\,Proposition 3.1 was inspired by [1, Proposition
1.1].}

\begin{proposition} There exists a natural exact sequence
$$\begin{aligned}
\displaystyle{\bigoplus_{y\in C_{0}}}& CH^{d-1}(X_{y})^{\e\prime}\ra
{\rm{Ker}}\!\left[\,CH^{d}(X)^{\e\prime}\ra
CH^{d}(X_{\eta})^{\e\prime}\e\right]\\
&\ra {\rm{Ker}}\!\left[\e
E\!\left(\e\xb/\e\cb\,\right)\ra\displaystyle{\bigoplus_{y\in
C_{0}}}\,\,\frac{CH^{d-1}\!\left(\e
\xb_{\overline{y}}\right)^{\g_{\ng
y}}}{\pi_{\overline{y}}^{*}\,CH^{d-1}(X_{y})}\,\right]\ra 0\,,
\end{aligned}
$$
where $E\!\left(\e\xb/\e\cb\,\right)$ is the group (3).
\end{proposition}

As regards the right-hand group in the exact sequence of the
proposition, the following holds. Let
$$
H^{d-1}\!\left(X_{\eta},\e{\mathcal
K}_{d}\right)^{\prime}={\rm{Im}}\ng\left[H^{d-1}\!\left(X_{\eta},
\e{\mathcal K}_{d}\right)\ra
H^{d-1}\!\left(\xb_{\overline{\eta}},\e{\mathcal
K}_{d}\right)^{\g}\e\right]
$$
and
$$
{\rm{Sal}}_{\e d}(X/C)=\ng\left\{\ng f\!\in\!
H^{d-1}\!\left(\xb_{\overline{\eta}},\e{\mathcal
K}_{d}\right)^{\g}\!\colon\!\forall\e y\in
C_{0},\,\overline{\delta}_{\e\overline{y}}\e(f)\in\pi_{\overline{y}}^{*}\,
CH^{d-1}(X_{y})\right\},
$$
where $\overline{\delta}$ and $\pi_{\overline{y}}^{*}$ are the maps
of diagram (2).

\begin{proposition} There exists a natural exact sequence
$$\begin{aligned}
0&\ra\frac{{\rm{Sal}}_{\e d}(X/C)}{\left(\overline{\jmath}^{\e *}
H^{d-1}\ng\left(\e\xb,\e{\mathcal K}_{d}\right)\right)^{\g}\!\cdot\!
H^{d-1}(X_{\eta},\e{\mathcal
K}_{d})^{\e\prime}}\\
&\ra{\rm{Ker}}\!\left[\e
E\!\left(\e\xb/\e\cb\,\right)\ra\displaystyle{\bigoplus_{y\in
C_{0}}}\,\,\frac{CH^{d-1}\!\left(\e
\xb_{\overline{y}}\right)^{\g_{\ng
y}}}{\pi_{\ng\overline{y}}^{*}\,CH^{d-1}(X_{y})}\,\right]\\
&\ra H^{1}\!\left(\g,\e \overline{\jmath}^{\e *} H^{d-1}\ng\left(\e
\xb,\e{\mathcal K}_{d}\right)\right).
\end{aligned}
$$
\end{proposition}
\begin{proof} This follows by applying the snake lemma to a diagram
of the form
\[
\xymatrix{0\ar[r]&A\ar[r]\ar[d]&B\ar[r]\ar[d] & B/A \ar[d]\ar[r]&0\\
0\ar[r]&\overline{A}^{\,\g}\ar[r]&\overline{B}^{\,\g}\ar[r] &
\left(\e\overline{A}/\overline{B}\e\right)^{\,\g}\ar[r]&H^{1}\!\left(\e\g,
\overline{A}\e\right)&\\
}
\]
with $\overline{A}=\overline{\jmath}^{\,*} H^{d-1}\!\left(\e
\xb,\e{\mathcal K}_{d}\right)$, $\overline{B}=H^{d-1}\!\left(\e
\xb_{\overline{\eta}},\e{\mathcal K}_{d}\right)$, etc.
\end{proof}

\section{Proof of the main theorem.}
Let $C$ and $A$ be as in Section 2, let $p$ be a prime number
different from the characteristic of $k$ and let $X$ be a
Severi-Brauer scheme over $C$ of relative dimension $p-1$.

\begin{lemma} There exists a $\g$-isomorphism
$$
\overline{\jmath}^{\,*} H^{d-1}\!\left(\e \xb,\e{\mathcal
K}_{d}\right)\simeq \kc.
$$
\end{lemma}
\begin{proof} Clearly,
$\overline{\jmath}^{\,*} H^{d-1}\!\left(\e \xb,\e{\mathcal
K}_{d}\right)$ is the kernel of the map
$$
\overline{\delta}\colon
H^{d-1}\!\left(\e\xb_{\overline{\eta}},\e{\mathcal
K}_{d}\right)\ra\displaystyle\bigoplus_{\overline{y}\e\mid
y}CH^{d-1}\!\left(\e\xb_{\overline{y}}\e\right)
$$
appearing in the exact sequence (1) over $\kb$. Now
$\xb_{\overline{\eta}}\simeq{\mathbb P}^{\e p-1}_{\overline{\eta}}$
and $\xb_{\overline{y}}\simeq{\mathbb P}^{\e p-1}_{\kb}$ for every
$\overline{y}$, whence we have $\g$-isomorphisms
$$
H^{d-1}\!\left(\e\xb_{\overline{\eta}},\e{\mathcal
K}_{d}\right)\simeq \kb(C)^{*}
$$
and
$$
CH^{d-1}\!\left(\e\xb_{\overline{y}}\e\right)\simeq\mathbb Z
$$
for each $\overline{y}$. Under these isomorphisms, the map
$\overline{\delta}$ above corresponds to the canonical map
$$\begin{aligned}
&\kb(C)^{*}\ra\bigoplus _{\overline{y}\e\mid y}\,\mathbb Z,\\
&f\mapsto({\rm{ord}}_{\,\overline{y}}(f))_{\e\overline{y}\e\mid y},
\end{aligned}
$$
which yields the lemma.
\end{proof}

\begin{teorema} For every
$d$ such that $2\leq d\leq p$, there exists a canonical isomorphism
$$
CH^{d}(X)^{\e\prime}\simeq CH_{0}(X/C).
$$
\end{teorema}
\begin{proof} By Lemma 4.1, Hilbert's Theorem 90 and Proposition
3.2, there exists a natural isomorphism
$$
{\rm{Ker}}\!\left[\e
E\!\left(\e\xb/\e\cb\,\right)\ra\displaystyle{\bigoplus_{y\in
C_{0}}}\,\,\frac{CH^{d-1}\!\left(\e
\xb_{\overline{y}}\right)^{\g_{\ng
y}}}{\pi_{\ng\overline{y}}^{*}\,CH^{d-1}(X_{y})}\,\right]\simeq
\frac{{\rm{Sal}}_{\e d}(X/C)}{k^{*}H^{d-1}(X_{\eta},\e{\mathcal
K}_{d})^{\e\prime}}.
$$
On the other hand, by [5,(8.7.2)],
$H^{d-1}\!\left(X_{\eta},\e{\mathcal K}_{d}\right)^{\prime}=
{\rm{Nrd}}\e A^{*}$ for every $d$ such that $2\leq d\leq p$ and, for
each $y\in C_{0}$,
$$
\pi_{\overline{y}}^{*}\, CH^{d-1}(X_{y})\simeq
\pi_{\overline{y}}^{*}\, CH^{p-1}(X_{y})\simeq(q_{\e
y})_{*}CH_{0}(X_{y})\quad(=\mathbb Z\,\text{ or }\,p\e\mathbb Z).
$$
The latter implies that ${\rm{Sal}}_{\e
d}(X/C)=k(C)^{*}_{{\rm{dn}}}$, whence
$$\begin{aligned}
{\rm{Sal}}_{\e d}(X/C)/k^{*}H^{d-1}\!\left(X_{\eta},\e{\mathcal
K}_{d}\right)^{\e\prime}&\simeq
k(C)^{*}_{{\rm{dn}}}/k^{*}\e{\rm{Nrd}}\e A^{*}\\
&\simeq CH_{0}(X/C).
\end{aligned}
$$
Finally, [loc.cit.] shows that the groups $CH^{d}(X_{\eta})$ and
$CH^{d-1}\!\left(X_{y}\right)$ ($y\in C_{0}$) are torsion free,
whence $CH^{d}(X_{\eta})^{\prime}$ and
$CH^{d-1}\!\left(X_{y}\right)^{\prime}$ vanish. The theorem now
follows from Proposition 3.1.\end{proof}

\begin{corollary} Let $d$ be such that $2\leq d\leq p$. Then
$CH^{d}(X)^{\e\prime}$ is finite if
\begin{enumerate}
\item $k$ is a number field, or
\item $k$ is a field of finite type over $\Bbb Q$, $C={\mathbb
P}^{1}_{k}$ and $X$ has a 0-cycle of degree one.
\end{enumerate}
In each of these cases, the group $CH^{d}(X)$ is finitely generated.
\end{corollary}
\begin{proof} Indeed, in these cases the group $CH_{0}(X/C)$
is finite [2].
\end{proof}


\begin{thebibliography}{7}



\bibitem[1]{1} Colliot-Th\'el\`ene, J.-L. and Skorobogatov, A.:\emph{
 Groupe de Chow des z\'ero-cycles sur les fibr\'es en
quadriques.} $K$-Theory {\bf{7}}, pp. 477-500 (1993).

\bibitem[2]{2} Frossard, E.:\emph{ Groupe de Chow de dimension
z\'ero des fibrations en vari\'et\'es de Severi-Brauer.} Comp. Math.
{\bf{110}}, pp. 187-213 (1998).

\bibitem[3]{3} Fulton, W.:\emph{ Intersection Theory.} Second Ed.
Springer-Verlag, 1998.


\bibitem[4]{4}  Grothendieck, A.:\emph{ Le Groupe de Brauer I.} In:
Dix Expos\'es sur la Cohomologie des Sch\'emas. North-Holland,
Amsterdam, pp. 46-66 (1968).

\bibitem[5]{5}  Merkurjev, A.S. and Suslin, A.A.:\emph{
$K$-cohomology of Severi-Brauer varieties and the norm residue
homomorphism.} Math. USSR Izv. {\bf{21}}, No. 2, pp. 307-340 (1982).



\bibitem[6]{6} Salberger, P.:\emph{ Galois descent and class groups
of orders.} Lect. Notes in Math. {\bf{1142}}, pp. 239-255 (1985).

\bibitem[7]{7} Sherman, C.:\emph{ Some theorems on the $K$-Theory
of coherent sheaves.} Comm. in Alg. {\bf{7}}(14), pp. 1489-1508
(1979).


\end{thebibliography}
\end{document}